\documentstyle[12pt]{article} 
\newfont{\twlvmsb}{msbm10 scaled\magstep1}
\newfont{\ninemsb}{msbm9}
\newfont{\sixmsb}{msbm6}
\newfam\msbfam
\textfont\msbfam=\twlvmsb
\scriptfont\msbfam=\ninemsb
\scriptscriptfont\msbfam=\sixmsb
\catcode`\@=11
\def\Bbb{\ifmmode\let\next\Bbb@\else
  \def\next{\errmessage{Use \string\Bbb\space only in math mode}}\fi\next}
\def\Bbb@#1{{\Bbb@@{#1}}}
\def\Bbb@@#1{\fam\msbfam#1}
\newfont{\largeeufm}{eufm10 scaled\magstep4}
\newfont{\twlveufm}{eufm10 scaled\magstep1}
\newfont{\elveufm}{eufm10 at 11pt}
\newfont{\teneufm}{eufm10}
\newfont{\nineeufm}{eufm9}
\newfam\eufam
\textfont\eufam=\twlveufm
\scriptfont\eufam=\nineeufm
\def\frak{\ifmmode\let\next\frak@\else
\def\next{\errmessage{Use \string\frak\space only in math mode}}\fi\next}
\def\frak@#1{{\fam\eufam{{#1}}}}
\catcode`@=12
\newcommand{\Z}{{\Bbb Z}} 
\newcommand{\Zz}{{\Bbb Z}_2 } 
\newcommand{\C}{{\Bbb C}} 
\newcommand{\LC}{{\bf\Lambda}_\infty} 
\newcommand{\gl}{gl(m|n)} 
 
\newcommand{\Uq}{U_q(gl(m|n))} 
\newcommand{\DE}{\overline{E}} 
 
\newcommand{\ba}{\begin{eqnarray}}
\newcommand{\na}{\end{eqnarray}}
\newcommand{\ban}{\begin{eqnarray*}}
\newcommand{\nan}{\end{eqnarray*}}
\newtheorem{lemma}{Lemma}
\newtheorem{proposition}{Proposition}

\newcommand{\BPW}{Poincar\'e - Birkhoff - Witt\ } 
\newcommand{\BBW}{Bott - Borel - Weil\ } 

\begin{document} 
\title{\normalsize{\bf  BOTT - BOREL - WEIL CONSTRUCTION FOR \\ 
QUANTUM SUPERGROUP $U_q(gl(m|n))$ } } 
\author{\small R. B. ZHANG \\
\small Department of Pure Mathematics,  University of Adelaide, 
Adelaide, Australia} 
\date{} 
\maketitle

\begin{abstract}  
The finite dimensional irreducible representations of the 
quantum supergroup $U_q(gl(m|n))$ are constructed geometrically 
using techniques from the Bott - Borel - Weil theory and vector 
coherent states. 
\end{abstract} 

\vspace{2cm}

\noindent
\section{\normalsize INTRODUCTION}\normalsize  
Supersymmetry and quantum groups are two of the most important 
discoveries of mathematical physics in recent times. These two 
notions have been combined together in the theory of quantum 
supergroups, which has been under intensive investigation in 
the last few years. The origin of quantum supergroups may be 
traced back to the Perk - Schultz model  in statistical mechanics, 
but the systematic study of these remarkable algebraic 
structures only began about six years ago. Since then much has been 
understood about  both the structures and representations of 
the quantum supergroups, and their applications have also 
been widely explored, leading to significant advances in various areas 
of physics and mathematics, notably, integrable models and low 
dimensional topology. 
Amongst the quantum supergroups arising from  `deformations' 
of the universal enveloping algebras of basic classical Lie 
superalgebras, $U_q(gl(m|n))$ has been  best studied, in particular, 
its representation theory was developed in \cite{II}. 
The present paper aims to further 
develop the representation theory of $U_q(gl(m|n))$ using the 
techniques from \BBW theory\cite{Bott}
and vector coherent states \cite{Hecht}. 

The \BBW theorem\cite{Bott} relates finite dimensional 
irreducible representations of compact Lie groups to 
cohomology groups of homogeneous vector bundles. 
Relevant to our investigations in this paper is the following 
situation: Let $G$ be a simple complex Lie group, which may be regarded 
as the complexification of some compact Lie group. 
Let $P\subset G$ be a parabolic subgroup. Given finite dimensional irrep 
$V_0$ of $P$ with a  highest weight which is assumed to be integral 
and dominant with respect to $G$ as well,  
one forms the homogeneous vector bundle 
$G\times _P V_0$ $\rightarrow$ $ G/P$.  
The holomorphic sections of the vector bundle 
furnishes a finite dimensional irrep of $G$. 
The method of vector coherent states, which was developed  
to address specific problems in quantum mechanics, 
provides a way to explicitly construct the holomorphic sections,
and to realize the Lie algebra of $G$ as differential operators 
on $G/P$,  thus putting the abstract \BBW construction 
in a concrete form.

A supersymmetric version of \BBW theorem  was investigated by 
Penkov and Serganova \cite{Penkov}\cite{Serganova}, who also 
extensively developed the theory of homogeneous superspaces.  
It is well known that supermanifold geometry is much 
richer than ordinary geometry. This makes the \BBW theory  
for Lie supergroups a very interesting, but also difficult, 
subject to study. 
 
In recent years, the \BBW theory has also been investigated for 
quantum groups \cite{Parshall}\cite{Lohe}\cite{Biedenharn}. 
In particular, the work of Biedenharn and Lohe \cite{Lohe}\cite{Biedenharn} 
on the quantized universal enveloping algebra $U_q(gl(m))$ is closely 
related to vector coherent states.  It has  
the virtue of being explicit and thus readily applicable 
to investigationts of mathematical structures relevant to physics, 
e.g.,  quantum group tensor operators\cite{Cornwell}.

In this paper we will carry out the \BBW construction for the quantum 
supergroup $\Uq$.  
As preparations for treating the quantum supergroup, 
and also for the purpose of understanding the 
underlying geometrical structure of vector coherent states,   
we study them extensively for classical $GL(m|n)$. 
In the quantum case, 
we realize the module of any finite dimensional irreducible 
representation of $U_q(gl(m|n))$ in terms of vector valued functions 
on some supermanifold, and the generators of the algebra by 
{\em difference operators} acting on these functions. 

The main body of the paper is divided into sections 2 and 3, 
respectively dealing with the \BBW construction for classical 
and quantum $GL(m|n)$. 
The first two subsections of section 2 are of a review nature, 
summarizing properties of structural  and representation 
theoretical features of the Lie superalgebra $\gl$. 
Some of the material covered can not be easily found in the 
mathematical physics literature, but  
is necessary for the remainder of the paper. 
Subsections 2.3 and 2.4 construct two types of vectors coherent states 
for the classical Lie supergroup $GL(m|n)$.  One kind is associated 
with the subgroup $GL(m|n-1)\times GL(1)$, and the other provides a 
new method for studying irreducible components of Kac modules. 
Subsections 2.5 elucidates the geometrical structure of the vector 
coherent states, thus to position the results of the previous two 
subsections into the general framework 
of \BBW theory for Lie supergroups. 

Section 3 consists of three subsections.  Subsection 3.1 investigates 
some structural and representation theoretical aspects of $\Uq$. 
Subsection 3.2 constructs the vector coherent states
for a special class of representations of the quantum supergroup, 
namely, the contragredient tensor irreps.  Some rather technical 
results needed for the next subsection are also proved here. 
The last subsection presents the \BBW construction 
of all the finite dimensional irreps of $\Uq$.

\section{\normalsize CLASSICAL  GL(m$\mid$n)}  
\subsection{ gl(m$\mid$n)}\normalsize 
Let us begin by introducing the familiar formulation of
Lie superalgebras given in \cite{Kac}\cite{Scheunert}.  
For the purpose of studying the representation theory 
of Lie supergroups and their quantum analogues in geometrical terms, 
one needs to define supergroups within the framework of 
supermanifold theory, and to formulate Lie superalgebras accordingly. 
However,  we will postpone this until
subsection 2.3 when we construct supersymmetric coherent states, 
where the theory of supermanifolds becomes indispensable.
At the algebraic level, the formulation of \cite{Kac}\cite{Scheunert} 
is much easier to handle.  

Within this formulation,  a Lie superalgebra is considered as a 
$\Zz$ graded Lie algebra over the complex field $\C$, 
namely, a $\Zz$ graded vector space endowed with a graded 
bracket.   The underlying vector space of the Lie superalgebra $\gl$ 
has the standard homogeneous basis  $\{ e_{a b} \ |\   a, b\in{\bf I}\}$, 
where the index set ${\bf I}$ is $\{1, 2, ..., m+n\}$.
Set ${\bf I}'={\bf I}\backslash\{m+n\}$. 
Introduce the gradation index $[\ ]:  {\bf I}\rightarrow
\Zz$ such that
$[a]=\left\{\begin{array}{l r}
          0,& a\le m,\\
          1,& a>m.
           \end{array} \right. $
Let $\gl_\eta$, $\eta\in \Zz$,  be the vector space over $\C$
spanned by the $e_{a b}$ with $[a]+[b]\equiv$ $\eta \ (mod\ 2)$.  
Then $\gl_0$ and $\gl_1$ are the even and odd subspaces of $\gl$ 
respectively. We will abuse the notation somewhat  
and define $[\ ]: \gl_0\cup \gl_1$  $\rightarrow$
$\Zz$, $[x]=\left\{\begin{array}{l r}
          0,& x\in\gl_0,\\
          1,& x\in\gl_1.
           \end{array} \right. $
Then the $\Zz$ graded Lie bracket for $\gl$ is defined by 
\ba
{[} e_{a b}, \ e_{c d}\} &=& \delta_{b c} e_{a d} -
(-1)^{([a]+[b])([c]+[d])}  e_{c b}\delta_{d a}.
\na
For convenience, we will regard $\gl$ as embedded in its universal
enveloping algebra. Thus the graded bracket $[\  , \ \}$ can be 
interpreted as the graded commutator
\ba
[x, \ y\}&=&x y - (-1)^{[x][y]} y x.\label{commutator}
\na

Let $\frak h$ be the Lie subalgebra generated by 
$\{ e_{a a}\, | \, a\in {\bf I}\}$, which is a  Cartan 
subalgebra of $\gl$. Introduce a basis $\{\epsilon_a \, | \, 
a\in {\bf I}\}$ for the dual vector space ${\frak h}^*$ such that 
$\epsilon_a ( e_{b b})=\delta_{a b}$. Note that the bilinear form 
for $\gl$ defined by  
\ban 
e_{a b} \otimes e_{c d}&\mapsto& (-1)^{[a]} \delta_{b c}\delta_{a d}, 
\nan 
is super invariant and nondegenerate. Its restriction to 
${\frak h}$ induces the following nondegenerate bilinear form  
\ban 
(\ , \ ): {\frak h}^*\otimes {\frak h}^* &\rightarrow& \C, \\
(\epsilon_a, \ \epsilon_b)& =& (-1)^{[a]}\delta_{a b}. 
\nan 

The positive and negative root spaces of $\gl$ with  respect to 
$\frak h$ are respectively given by
\ban 
{\frak n}^+&=& \bigoplus_{a< b} \C e_{a b}, \\ 
{\frak n}^-&=& \bigoplus_{a>b} \C e_{a b},
\nan   
which are nilpotent super subalgebras of $\gl$. We will denote by 
${\frak b}^{\pm}$ the Borel subalgebras $\frak h 
+ {\frak n}^\pm$ respectively.  
 
A parabolic subalgebra of $\gl$ is a proper subalgebra 
containing a Borel subalgebra. We call a parabolic subalgebra 
standard if it contains ${\frak b}^+$.  A feature distinct 
from those of ordinary Lie groups is that 
not all the parabolic subalgebras are Weyl group conjugate 
to standard ones. 

Now every standard parabolic subalgebra is of the form 
\ban 
{\frak p} &=& {\frak b}^+ + {\frak u}_-, 
\nan 
where ${\frak u}_-$ is a subalgebra of ${\frak n}^-$ generated by 
the elments of $\{ e_{ a+1\, a} \ \mid \  a \in\bf\Theta\}$, 
where $\bf\Theta$ is a proper subset of ${\bf I}'$.   
Note that $\bar{\frak u}_-={\frak n}^-\backslash {\frak u}_-$ 
is also a subalgebra of ${\frak n}^-$. 
For any Lie super subalgebra $\frak k$ of $\gl$, 
we denote by $U(\frak k)$ its universal enveloping algebra.  
Then it follows from the \BPW theorem that 
\ban 
U(\gl)&=& U(\bar{\frak u}_-) U(\frak p). 
\nan

\subsection{Irreducible representations}\normalsize 
Let $\frak p$ be any standard parabolic subalgebra, and 
$V_0(\lambda)$ be a finite dimensional irreducible $\frak p$ module. 
Since ${\frak h}$ is an abelian Lie subalgebra of $\frak p$, 
Lie's theorem asserts that there exists at least one common 
eigenvector of all elements of ${\frak h}$ in $V_0$. 
By noting the fact that $U(\frak p)$ can be decomposed into 
eigenspaces of ${\frak h}$ ( under the adjoint action ), we immediately 
see that the irreducible $\frak p$ module $V_0(\lambda)$ must be of highest 
weight type, i.e., there exists a  non - null 
$v_+\in V_0(\lambda)$ such that 
\ban 
e_{a\, b}\  v_+ &=&0, \ \ \ e_{a\, b}\in{\frak p}, \ a< b, \\ 
e_{a\, a} v_+ &=&\lambda_a v_+,   \ \ \ e_{a\, a}\in{\frak p},    
\nan 
and $V_0(\lambda)$ is cyclically generated by $v_+$.  

Construct the vector space 
\ban 
{\bar V}(\lambda)&=&U(\gl)\otimes_{U(\frak p)} V_0(\lambda).  
\nan 
Then clearly ${\bar V}(\lambda)=U(\bar{\frak u}_-)
\otimes V_0(\lambda)$,  and 
${\bar V}(\lambda)$ furnishes a $U(\gl)$ module with the action of 
$U(\gl)$ defined in the standard way:  for any $x\in U(\gl)$, 
$y\in U(\bar{\frak u}_-)$, $ x y$ can be expressed in the form 
\ban 
x y &=&\sum_t y_t' x_t',  \ \ \ x_t'\in U(\frak p), 
\ \ \ y_t'\in U(\bar{\frak u}_-).   
\nan  
Then
\ban  x\circ(y\otimes v_0)&=& \sum_t y_t'\otimes x_t' v_0,
 \ \ \ v_0\in V_0(\lambda). \nan   
${\bar V}(\lambda)$ can be decomposed into a direct sum of 
weight spaces, i.e., eigenspaces of $\frak h$, 
\ba 
{\bar V}(\lambda)&=& \bigoplus_{\omega\prec \lambda} {\bar V}^\omega, 
\ \ \ dim {\bar V}^\omega< \infty, \label{decomposition}  
\na 
where $\omega\in{\frak h}^*$, and $h v = \omega (h) v$, $\forall 
v\in {\bar V}^\omega$, $h\in \frak h$. The $\prec$ is a partial 
ordering of the elements of ${\frak h}^*$ defined by 
$\mu \prec\nu$ if $\mu - \nu$ $\in$ 
$\bigoplus_{a<b}\Z_+ (\epsilon_a -\epsilon_b)$.  
 
If $W$ is a proper submodule of ${\bar V}(\lambda)$, then 
$W\cap V_0(\lambda)=0$, 
as every nonvanishing element of $V_0(\lambda)$ cyclically generates  
${\bar V}(\lambda)$. Let $M$ be the union of all the proper submodules 
of ${\bar V}(\lambda)$.  
Then $M$ is again a proper submodule, which is 
unique and  is maximal in the sense that every other proper submodule 
of ${\bar V}(\lambda)$ is a submodule of $M$.  We will set $M=0$ 
in the case that no proper submodule of ${\bar V}(\lambda)$ exists. 
More generally, we define 
    \ban V(\lambda)&=& {\bar V}(\lambda)/M. \nan
It is a consequence of the maximality of $M$ that the quotient module 
$V(\lambda)$  is irreducible. Furthermore, $V(\lambda)$ 
admits a weight space decomposition similar to (\ref{decomposition}). 
Also, the canonical projection restricted to 
$V_0(\lambda)$ is one - to - one, and the image of $v_+$ is the 
unique maximal vector of $V(\lambda)$. 
      
The irreducible $U(\gl)$ module is finite dimensional if and only if  
\ba 
\lambda_a - \lambda_{a+1}\in \Z_+, \ \ \  m\ne a\in {\bf I}'.
\label{finite}  
\na 
This is the familiar finite dimensionality condition obtained by 
Kac. To understand it within our framework, we use the \BPW theorem 
for $U(\gl)$ again but in a different form to express it as 
\ban 
U(\gl)&=&U({\frak f}_-) U(gl(m)\oplus gl(n)) U({\frak f}_+), 
\nan
where $gl(m)\oplus gl(n)$ $\subset \gl$ is the maximal even 
subalgebra,  ${\frak f}_+$ is the subalgebra spanned by all the 
odd raising elements $e_{i \mu}$, $i\le m$, $\mu>m$, 
and ${\frak f}_-$ is that spanned by $e_{\mu i}$, $i\le m$, $\mu>m$. 
Note that both $U({\frak f}_\pm)$ are isomorphic to the grassmann 
algebra with $m n$ generators, hence $dim U({\frak f}_\pm)=2^{m n}$.  
Let $v^\lambda_+$ be the maximal vector of $V(\lambda)$.  
Consider $W_0=U(gl(m)\oplus gl(n)) v^\lambda_+$. 
As a $gl(m)\oplus gl(n)$ module, $W_0$ must be irreducible. 
Otherwise, a proper $gl(m)\oplus gl(n)$ submodule $W_0'$ of 
$W_0$ would generate a proper $\gl$ submodule $U({\frak f}_-)W_0'$ 
$\subset$  $V(\lambda)$, contradicting the irreducibility of $V(\lambda)$. 
Now $W_0$ is finite dimensional if and only if 
$\lambda$ satisfies (\ref{finite}), and this in turn leads to our claim. 

Recall that every finite dimensional irrep of $\gl$ is of highest weight 
type, and is uniquely determined by a highest weight 
$\lambda\in{\frak h}^*$ satisfying (\ref{finite}). 
Therefore, for any given standard parabolic subalgebra $\frak p$, 
the construction presented above yields all the finite dimensional 
irreps of $\gl$. 
It is worth noting that when $\frak p = \gl\backslash {\frak f}_-$, 
the construction coincides with that of Kac.

\subsection{Vector coherent states: 1} 
{}From this subsection on, we need to formulate Lie superalgebras 
over supernumbers, following  the treatment of DeWitt \cite{DeWitt}. 
Introduce the infinite complex Grassmann 
algebra ${\bf\Lambda}_\infty$ completed with Frechet topology. 
We will denote by $\C_c$ the commutative part, and by 
$\C_a$ the anti - commutative part of ${\bf\Lambda}_\infty$,  
and  represent $\C_a^m\times \C_c^n$ by $\C^{m|n}$.    

We can now define a Lie superalgebra as a supervector space over 
${\bf\Lambda}_\infty$ endowed with a Lie super bracket. In particular, 
the Lie super bracket of $gl(m|n)$ is the same as that given in 
subsection 2.1. However, it is worth stressing that even elements 
of $gl(m|n)$ commute with all supernumbers, while odd elements 
anti - commute with $\C_a$. 
Our discussions on structures and representations of $gl(m|n)$ 
in the previous two subsections are still valid in this more 
general setting. 

Let $\frak p$ be the parabolic subalgebra 
spanned by the following elements
\ban 
e_{a b},& \ a, b\in{\bf I},& a\le b,\\ 
e_{d c}, & c, d\in{\bf I}', & \ c<d.
\nan 
Note that the $gl(m| n-1)\oplus gl(1)$ subalgebra spanned 
by $\{e_{a b} \ | \ a, b\in {\bf I}'\}$ 
$\cup \{e_{m+n\, m+n}\}$ is contained in ${\frak p}$. 
Let $V_0(\lambda)$ be a finite dimensional irreducible $\frak p$ module 
with highest weight $\lambda\in{\frak h}^*$ satisfying the condition 
(\ref{finite}), and denote by $\pi_0:  {\frak p}\rightarrow$ 
$gl(V_0(\lambda))$ the associated irrep of ${\frak p}$. 
The induced module construction presented in the 
last subsection allows us to construct a finite dimensional irreducible 
$\gl$ module $V(\lambda)$ into which $V_0(\lambda)$ is canonically 
embedded.  We can decompose $V(\lambda)$ into a direct sum 
of $e_{m+n\, m+n}$ eigenspaces 
\ban 
V(\lambda)&=&\bigoplus_{k=0}^L V_k, 
\nan 
where $e_{m+n\, m+n}$ takes eigenvalue $\lambda_{m+n} + k$ 
when acting on $V_k$. Note that $V_0=V_0(\lambda)$.

Let $\{ w^\alpha \ | \ \alpha = 1, 2, ..., d\}$ 
be a basis of $V(\lambda)$ such that $\{ v^i=w^i   \ | \ i=1, 2, ..., 
d_0\}$ forms a basis of $V_0(\lambda)$, where $d= dim V(\lambda)$ , 
and $d_0=dim V_0(\lambda)$.  
Introduce the basis $\{ \langle w^\alpha\mid  \ 
| \ \alpha = 1, 2, ..., d\}$  for the dual vector 
space $V(\lambda)^*$ of $V(\lambda)$ such that 
\ban 
\langle w^\alpha \mid w^\beta\rangle &=& \delta_{\alpha \beta}, 
\nan 
where $\langle\  \mid\  \rangle$ denotes the dual space pairing. 
Then the matrix elements of any $X\in\gl$ in the irrep 
$\pi: \gl\rightarrow gl(V(\lambda))$ furnished by the module 
$V(\lambda)$ are given by 
\ban 
\pi(X)_{\alpha  \beta} &=& \langle w^\alpha \mid X w^\beta\rangle. 
\nan  
It is a standard result that $V(\lambda)^*$ carries a left $\gl$ module 
structure defined by 
\ban 
(X\circ \langle w_1\mid)w_2 \rangle &=& -\langle w_1\mid X w_2 \rangle.  
\nan 
Introduce $( Z_a )_{a\in{\bf I}'}$ $\in \C^{m|n-1}$. Set 
$Z_i=\theta_i$, $1\le i\le m$, $Z_\mu=z_\mu$, $m<\mu<m+n$. 
Then the $z_\mu$ are even,   while the $\theta_i$ are odd, 
and obey $\theta_i \theta_j = - \theta_j  \theta_i$. 
Denote by $\LC[[Z]]$ the ring of polynomials in the variables 
$Z_a$, and $\LC[[Z]]_L$ the subset of the  polynomials of degree 
less or equal to $L$, which is clearly finite dimensional.  
Needless to say, it is assumed that $z_\mu$ commutes with all 
elements of $U(\gl)$, while $\theta_i$ commutes with the even 
elements and  anticommutes with the odd ones.   
We construct the linear space $V(\lambda)[[Z]]$ $=$ 
$\LC[[Z]]\otimes V(\lambda)$, and define the bilinear map 
$V(\lambda)^* \otimes V(\lambda)[[Z]]$ $\rightarrow$ 
$ \LC[[Z]]$  by 
\ban 
\langle w_1 | p(Z)\otimes w_2\rangle &=& 
(-1)^{[w_1][p(Z)]} \langle w_1|w_2\rangle p(Z). 
\nan  
Set 
\ban g(Z)&=& \exp{\left( \sum _{a\in{\bf I}'} Z_a 
     \otimes e_{a m+n}\right)}, \nan 
where the exponential should be understood as a formal power series 
at this stage. Now we define a linear map $\xi:$  $V(\lambda)$ 
$\rightarrow$   $\LC[[Z]]_L\otimes V_0(\lambda)$ by 
\ba  
\xi_w(Z)&=& \sum_{i=1}^{d_0} (-1)^{[v_i](1+[w])}
       \langle v^i |  g(Z) ( 1 \otimes w) \rangle \otimes v^i. 
\na    
We set $\xi_\alpha(Z)=\xi_{w^\alpha}(Z)$.  
Denote the supervector space spanned by
all the $\xi_\alpha(Z)$ by $V^\lambda_Z$, which will be called the 
space of vector coherent states. Then 
\begin{proposition} 
\begin{enumerate} 
\item The $\xi_\alpha(Z)$, $\alpha = 1, 2, ..., d$, are linearly 
independent over $\LC$. 
\item $V^\lambda_Z$ furnishes a $\gl$ module with the action defined by 
\ban 
X\circ \xi_w(Z) &=& \sum_{i=1}^{d_0} (-1)^{[v_i](1+[X]+[w])}
       \langle v^i |  g(Z) ( 1 \otimes X w) \rangle \otimes v^i, \\ 
\mbox{or more compactly,} \\ 
X\circ \xi_w(Z) &=& \xi_{X w}(Z), \ \ \   \forall X\in \gl. 
\nan 
\item $V^\lambda_Z $ and  $ V(\lambda)$ are isomorphic as $\gl$ modules.    
\end{enumerate}
\end{proposition} 
{\em Proof}:  Part 1 is an immediate consequence of the irreducibility 
of $V(\lambda)$, while part 2 follows from part 1 and the following simple 
calculation
\ban 
X\circ \xi_\alpha(Z)&=& \sum_{\beta} \pi(X)_{\beta \alpha} \xi_\beta(Z). 
\nan 
Observe that $1\otimes v^\lambda_+$ is the maximal vector of 
$V^\lambda_Z $.  It is easy to see that the weight of $1\otimes v^\lambda_+$ 
is $\lambda$, hence our claim in part 3.

Now we have realized a finite dimensional irreducible $\Uq$ module as 
a subset of $\LC[[Z]]_L\otimes V_0(\lambda)$.  Under a suitable 
dualization,  elements of $V^\lambda_Z $ may be interpreted as 
holomorphic sections of a homogeneous supervector bundle. 
We will give more details on this later. 
Here we present the explicit realization of $\gl$ 
in terms of differential operators on the superbundle.
\begin{proposition} 
The irreducible representation $\pi: \gl\rightarrow End(V^\lambda_Z)$ of
$\gl$ can be realized by 
\ba
\pi(e_{a\,  b}) &=& -(-1)^{[a]([b]+1)} Z_b {{\partial}\over{\partial Z_a}} 
\otimes 1 + 1\otimes \pi_0(e_{a\, b}),\nonumber \\ 
\pi(e_{m+n\, m+n})&=& - \sum_{a\in{\bf I}'}Z_a 
{{\partial}\over{\partial Z_a}}\otimes 1 + 1\otimes \pi_0(e_{m+n\, m+n}), 
\nonumber \\ 
\pi(e_{m+n\, a})&=& \sum_{b\in{\bf I}'} 
              \left\{Z_b\otimes \pi_0(e_{b a}) 
              + (-1)^{[a]}Z_a Z_b {{\partial}\over{\partial Z_b}}
              \otimes 1 \right\} 
              + Z_a\otimes \pi_0(e_{m+n\, m+n}),  \nonumber \\ 
\pi(e_{a \ m+n})&=& {{\partial}\over{\partial Z_a}}\otimes 1, 
\hspace{3cm}  \ a, b\in{\bf I}'.     \label{differential}  
\na 
\end{proposition}  
{\em Proof}:  The correctness of the realization for $e_{a \, m+n}$ is 
clear. For the $e_{a\, b}$, $a, b\in{\bf I}'$, we observe that 
\ban -(-1)^{[a]([b]+1)} Z_b {{\partial}\over{\partial Z_a}}
\otimes 1 + 1\otimes e_{a\, b}\nan       
commutes with $g(Z)$. Hence 
\ban 
e_{a\, b}\circ \xi_\alpha(Z) &=&\sum_{i=1}^{d_0} 
\langle v^i |\left\{-(-1)^{[a]([b]+1)} 
Z_b {{\partial}\over{\partial Z_a}}
\otimes 1 + 1\otimes e_{a\,  b}\right\}   
g(Z) ( 1 \otimes  w) \rangle \\ 
{} & \otimes & v^i (-1)^{[v_i]([a]+[b]+1+[w])}, 
\nan 
which immediately yields what we want to prove.  
The case of $e_{m+n\, m+n}$ 
can be shown similarly. As to $e_{m+n\, a}$, by using 
\ban 
g(Z)(1\otimes e_{m+n\, a})g(Z)^{-1}&=& 1\otimes e_{m+n\, a} 
+ (-1)^{[a]} Z_a \sum_{b\in{\bf I}'} 
Z_b {{\partial}\over{\partial Z_b}} \otimes 1  \\
&+& (-1)^{[a]} Z_a\otimes e_{m+n\, m+n} 
+ \sum_{b\in{\bf I}'} Z_b\otimes e_{b\, a},   
\nan 
we obtain   
\ban 
 e_{m+n\, a}&\circ &\xi_\alpha(Z)= 
\left\{(-1)^{[a]} Z_a \sum_{b\in{\bf I}'} 
Z_b {{\partial}\over{\partial Z_b}} \otimes 1 \right\}\xi_\alpha(Z)\\ 
&+& \sum_{i,  j} \langle v^j |\left\{ 
(-1)^{[a]} Z_a\pi_0(e_{m+n\, m+n})_{i j}       
+ \sum_{b\in{\bf I}'} Z_b\pi_0(e_{b a})_{i j}\right\} 
g(Z)(1\otimes w^\alpha)\rangle\\ 
&\otimes& v^i (-1)^{[v_i]([a] + [w])}, 
\nan 
from which the last formula of (\ref{differential}) can be read off.  
 
Several features of this realization  are worth observing. 
Repeated applications of $\pi(e_{a b})$, $a<b$, $a, b\in{\bf I}$, to  
$1\otimes v^\lambda_+$, the highest weight vector of $V^\lambda_Z$, 
generate the irreducible module $V^\lambda_Z$ automatically, rather 
than an analogue of ${\bar V}(\lambda)$.  Also, the problem of 
studying the irreps of $\gl$ may now be dealt with by studying the 
irreps of the chain of subalgebras  
\ban 
\gl\supset gl(m|n-1)\supset ... \supset gl(m|1)\supset gl(m). 
\nan 
This may provide a practical method for investigating the
detailed structures of the finite dimensional irreps of $\gl$.

\subsection{Vector coherent states: 2}\normalsize 
As we have already pointed out, Kac' induced module construction 
is a special case of the method developed in subsection 2, with 
the parabolic subalgebra chosen to be
\ban 
{\frak p} &=& gl(m) + gl(n) + {\frak f}_+. 
\nan 
The vector coherent states associated with this parabolic subalgebra 
were studied by Le Blanc and Rowe\cite{Rowe}. Here we briefly 
outline the construction. 

Let $V_0(\lambda)$ be a finite dimensional irreducible 
$\frak p$ module with highest weight $\lambda$. Following the
general procedure of subsection 2, we construct ${\bar V}(\lambda)$, 
which is isomorphic to $U({\frak f}_-)\otimes V_0(\lambda)$, 
thus is a Kac module.  As before, we denote by $M$ the 
maximal proper submodule of ${\bar V}(\lambda)$ in the atypical 
case, and set $M=0$ in the typical case. 

Let $( \theta_{\mu i} )_{1\le 1\le m,\, m+1\le \mu\le m+n}$ be 
any point in $\C^{m n|0}$. Define  
\ban 
g(\theta)&=& exp{\left(\sum_{i=1}^m\sum_{\mu=m+1}^{m+n} 
\theta_{\mu i}\otimes e_{i \mu}\right)}. 
\nan 
We define a linear map $\Psi: {\bar V}(\lambda)\rightarrow 
\LC[[\theta]]\otimes V_0(\lambda)$ by   
\ban 
w\mapsto \xi_w(\theta)&=&\sum_{i=1}^{d_0} (-1)^{[v_i](1+[w])} \langle v^i 
| g(\theta) ( 1 \otimes w) \rangle \otimes v^i.
\nan 
Then $Im\Psi$ furnishes a $\gl$ module with the module action 
\ban 
X\circ \xi_w(\theta)&=&\sum_{i=1}^{d_0} (-1)^{[v_i](1+[X]+[w])} \langle v^i
| g(\theta) ( 1 \otimes X w) \rangle \otimes v^i,   \ \ \ \forall X\in\gl. 
\nan 
Furthermore,  
\begin{proposition}
\begin{enumerate} 
\item $Ker\Psi=M$, and $Im\Psi$ is irreducible. 
\item When acting on $Im\Psi$, $\gl$ is realized by 
\ba
\pi(e_{i j}) &=& -\sum_{\mu=m+1}^{m+n} \theta_{j \mu} 
{{\partial}\over{\partial \theta_{i \mu} } } \otimes 1 
+ 1\otimes \pi_0(e_{i j}), \nonumber \\
\pi(e_{\mu \nu}) &=& -\sum_{i=1}^{m} \theta_{i \nu} 
{{\partial}\over{\partial \theta_{i \mu}}} \otimes 1 
+ 1\otimes \pi_0(e_{\mu \nu}),  \nonumber \\
\pi(e_{\mu i})&=& \sum_{j=1}^{m}\sum_{\nu=m+1}^{m+n} 
   \theta_{\nu i}\theta_{\mu j} 
{{\partial}\over{\partial \theta_{\nu j} }}\otimes 1\nonumber\\
&+& \sum_{j=1}^{m}\theta_{\mu j}\otimes\pi_0(e_{j i})
+\sum_{\nu=m+1}^{m+n}\theta_{\nu i}\otimes\pi_0(e_{\mu \nu}), 
\nonumber \\  
\pi(e_{i \mu})&=& {{\partial}\over{\partial \theta_{\mu i}}}\otimes 1,
\hspace{3cm}  i, j\le m, \ \ \ \mu, \nu>m,     
\na 
where $\pi_0: gl(m)\oplus gl(n)\rightarrow gl(V_0(\lambda))$ is the 
irrep of $gl(m)\oplus gl(n)$ afforded by $V_0(\lambda)$. 
\end{enumerate}
\end{proposition}  
Results of this Proposition for the special case of $gl(m|1)$ 
were known to Dundi and Jarvis \cite{Jarvis}, and the general case 
were contained in some unpublished work of Bracken's \cite{Bracken}, 
and the work of Le Blanc and Rowe\cite{Rowe}.

In the next subsection we will demonstrate that the  
underlying geometry of the vector coherent states of this subsection, 
and thus also Kac' induced module construction, is very special.

\subsection{Geometrical interpretation} 
The geometrical meaning of the vector coherent states for $gl(m|n)$ 
will be elucidated in this subsection. For this purpose, we need to  
consider supermanifolds more general than $\C^{m|n}$. We will 
follow DeWitt's geometrical formulation,  within which
 a supermanifold is, roughly speaking,  
a collection of open sets of $\C^{m|n}$ patched together, where 
a set is open if it is the pre - image of an open set in 
$\C^n$ under the natural projection of $\C^{m|n}$ onto $\C^n$.

As is well known, there are several distinct formulations 
of supermanifolds, but the most widely studied ones are those 
due respectively to DeWitt and to Berezin - Kostant - Leites. 
While  DeWitt's is suitable for theoretical physics,  
the sheaf theoretical approach of  
Berezin - Kostant - Leites has many mathematically appealing 
features.  It has also been shown  \cite{Bartocci} that 
the DeWitt and Berezin - Kostant - Leites supermanifolds are 
the same in an appropriate sense. 
We may also mention that we will only deal with rather 
simple supermanifolds, and do not concern ourselves with 
detailed functional analysis on them. 
Therefore the various approaches probably do 
not lead to any essential differences for us here.  
 
Let $G$ be the component of the Lie supergroup $GL(m|n)$ connected 
to the identity. Needless to say, $G$ is  a  Lie supergroup in its own 
right, and has the Lie superalgebra $\gl$.  
Let $P\subset G$ be the Lie super subgroup with Lie superalgebra 
$\frak p$, which is assumed to be a parabolic subalgebra of $\gl$. 

We first consider the case where the parabolic subalgebra ${\frak p}$ 
is generated by $\{ e_{a b},\ e_{m+n\, a},$  $\ e_{m+n\, m+n},$ 
$a, b\in{\bf I}'\}$. Note that ${\frak p}\supset {\frak b}^{-}$. 
Now $G/P$ (understood as the right coset space. 
A more appropriate notation may be $P\backslash G$.) yields a 
homogeneous superspace \cite{Manin} , 
which coincides with the projective superspace 
$\C P^{m|n-1}$ defined in the following way.  
Let $( W_a )_{a\in{\bf I}}$ $=$ $(\zeta_i; w_\mu)_{1\le i\le m;
m+1\le \mu\le m+n}$ be any point of the flat $m|n$ superspace $\C^{m|n}$. 
Consider the subspace $\C^{m|n}_*$ of $\C^{m|n}$ consisting of the points  
such that the bodies of the commuting coordinates of any 
given point do not vanish simultaneously.  Define an equivalence 
relation $( W_a )\sim ( W'_a )$,  $( W_a ), ( W'_a )\in \C^{m|n}_*$,   
if $( W_a )= c ( W'_a )$ for some $c\in\C_c$ with nonzero body.   
Then $\C P^{m|n-1}=$ $\C^{m|n}_*/\sim$.   
This supermanifold can be covered by the following 
coordinate patches 
\ban 
U_\mu&=&\{ (\theta^\mu_1, ..., \theta^\mu_m; z^\mu_{m+1}, ..., z^\mu_{\mu-1}, 
z^\mu_{\mu+1}, ..., z^\mu_{m+n} )\},  
\ \ \ \mu=m+1, ..., m+n,  \\  
\theta^\mu_i&=& \zeta_i/w_{\mu}, \ \ \ z^\mu_\nu=w_\nu/w_{\mu}, 
\ \ \ \mbox{the body of}\,  w_{\mu}\ne 0.  
\nan 
The transition functions can also be easily obtained.   

Consider a finite dimensional irreducible left $G$ module 
$V(\lambda)$ with highest weight $\lambda$, and denote the 
associated irrep by $T_\lambda$. 
Let $V(\lambda)^*$ be the dual $G$ module, and 
define $V_0(\lambda)^*\subset V(\lambda)^*$ to be the unique irreducible 
left $P$ submodule containing the lowest weight vector of $V(\lambda)^*$. 
Given a basis $\{ \langle v^i | \}$ of $V_0(\lambda)^*$, we express the 
left action of $P$ by 
\ban 
p\circ \langle v_i| &=&\sum_j L(p)_{j i} \langle v_j|, 
\ \ \ p\in{P}. 
\nan 
There is also a natural right $P$ module structure 
on $V_0(\lambda)^*$ with  
\ban 
\langle v_i| \circ p =\sum_j \langle v_j| R(p)_{i j}, 
& R(p)_{i j} = L(p)_{i j} ,  \ \ p\in P.  
\nan 
Let us now define the following $V_0(\lambda)^*$ valued 
functions on $G$
\ba 
\eta_w(g)&=& \sum_{i=1}^{d_0} (-1)^{[v_i](1+[w])}
       \langle v^i |  T_\lambda(g) w \rangle \otimes \langle v^i |, 
\ \ \ g\in G, \ w\in V(\lambda),   
\na
and denote by ${\cal O}_\lambda(G/P)$ their linear span.  
The right translation of $G$ defines 
a module action on  ${\cal O}_\lambda(G/P)$
\ban (g\circ \eta_w)(h) &=& \eta_w(h g),     \nan  
and the associated irrep of $G$ is isomorphic to $T_\lambda$.

A critical property of ${\cal O}_\lambda(G/P)$ to be observed is that  
\ba 
\eta_w(p g)&=& \eta_w(g) R(p^{-1}), \ \ \ \forall p\in P. \label{section}  
\na 
This allows us to interpret elements of ${\cal O}_\lambda(G/P)$ 
as global sections of the supervector bundle
\ban
G\times_{P} V^*_0(\lambda)\rightarrow G/P,
\nan
which is the quotient space $G\times V^*_0(\lambda)/\sim$, with the
equivalence relation defined by
\ban
(p g,  \eta)&\sim& (g, \ \eta R(p^{-1}) ),   \ \ \ p\in P.
\nan
 
Let ${\cal U}\subset G$ be a neighbourhood of the identity $e\in G$. 
We consider ${\cal O}_\lambda(G/P)$ restricted to ${\cal U}$.  
Differentiating the associated irrep of $G$ yields 
an irrep of the Lie superalgebra $\gl$, 
which is regarded as the left invariant vector fields on $G$: 
\ban 
(X\circ\eta_w)(g)&=& \sum_{i=1}^{d_0} (-1)^{[v_i](1+[X]+[w])}
     \langle v^i |  T_\lambda(g) X\circ w \rangle \otimes \langle v^i |, \\ 
& & X\in\gl, \  \  g\in{\cal U}.  
\nan
Recall the following decomposition of the Lie supergroup $G$,  
\ban 
G&=&P Q,  \\ 
Q&=&\left\{   \left( \begin{array}{l l}  
                 I& u\\
                 0& 1
                 \end{array}\right) 
      \right\}.   
\nan 
For each $g\in{\cal U}\subset G$, write $g=p q$, $p\in P\cap{\cal U}$, 
$q\in Q\cap{\cal U}$.   Then $\eta_w(g)=\eta_w(q)R(p^{-1})$.  
Also, every elements of $Q\cap{\cal U}$ can be expressed as  
$exp{( \sum _{a\in{\bf I}'} Z_a e_{a m+n})}$,  $( Z_a )\in \C P^{m|n-1}$. 
Therefore,  on $\cal U$, the vector space ${\cal O}_\lambda(G/P)$ 
is isomorphic to $V^\lambda_Z$, with the isomorphism given by the 
duality between $V_0(\lambda)$ and $V_0(\lambda)^*$. 
By comparing the highest weights, we can see that this also 
defines a $\gl$ module isomorphism.

The result of subsection 2.4. can be interpreted similarly, 
by choosing the parabolic subalgebra  ${\frak p}$ to be 
$gl(m)+gl(n)+{\frak f}_-$.  We will not repeat the details 
here, but merely point out that in this case, the homogeneous 
superspace $G/P$ is the flat superspace  $\C^{m n| 0}$, which 
does not have a body.

\section{\normalsize QUANTUM GL(m$\mid$n)} 
\subsection{$\Uq$ and its irreps} 
We will consider Jimbo's version of the quantum supergroup 
$U_q(gl(m\mid n))$ over ${\bf\Lambda}_\infty$.  
Fix $q\in\C_c$, {\em the body of which is assumed to 
be nonzero and not a root of unity}.  
$\Uq$ is a ${\Bbb Z}_2$-graded unital algebra  generated by
$\{K_a,  \ K_a^{-1}, \ a\in {\bf I};   \ 
E_{b\   {b+1}},$ $ \ E_{b+1,   b}, \ b\in {\bf I}'\}$,  
subject to the following relations  
\ba 
K_a K_a^{-1}=1, 
& & K_a^{\pm  1} K_b^{\pm 1} = K_b^{\pm 1}  K_a^{\pm 1}, \nonumber \\ 
K_a E_{b\ b\pm 1} K_a^{-1} &=& 
q_a^{\delta_{a b} -\delta_{a\  b\pm 1}} E_{b\  b\pm 1}, \nonumber \\ 
{}[E_{a\,  a+1},\,E_{b+1\,  b}\}& =& \delta_{a b}
(K_a K_{a+1}^{-1} - K_a^{-1} K_{a+1})/(q_a - q_a^{-1}),\nonumber \\ 
(E_{m\, m+1})^2 &=& (E_{m+1\, m})^2 = 0, \nonumber \\  
E_{a\,  a+1} E_{b\,  b+1} &=& E_{b\,  b+1} E_{a\,  a+1},\nonumber \\     
E_{a+1\, a} E_{b+1\, b} &=&E_{b+1\, b} E_{a+1\, a}, \ \ \  
\vert a - b\vert \ge 2, \nonumber \\ 
{\cal S}^{(+)}_{a \ a\pm 1}&=&{\cal S}^{(-)}_{a \ a\pm 1}=0,  
\ \ \ a\ne m,\nonumber \\ 
\{ E_{m-1\, m+2},\ E_{m\, m+1}\} &=& 
\{ E_{m+2\, m-1},\ E_{m+1\, m}\} = 0,  \label{quantum}  
\na 
where $q_a=q^{ (-1)^{[a]} }$,  
\ban 
{\cal S}^{(+)}_{a \ a\pm 1}&=& 
(E_{a\, a+1})^2  E_{a\pm 1\, a+1\pm 1} - (q +
q^{-1}) E_{a\, a+1} \ E_{a\pm 1\, a+1\pm 1} \ E_{a\, a+1}\\ 
& +& E_{a\pm 1\, a+1\pm1 }\ (E_{a\, a+1})^2,    \\ 
{\cal S}^{(-)}_{a \ a\pm 1}&=& 
(E_{a+1\, a})^2\,E_{a+1\pm 1\, a\pm 1} - (q +
q^{-1}) E_{a+1\, a}\ E_{a+1\pm 1\, a\pm 1} \ E_{a+1\, a}\\ 
&+& E_{a+1\pm 1\, a\pm 1}\ (E_{a+1\, a})^2, 
\nan 
and $E_{m-1\, m+2}$ and $E_{m+2\, m-1}$ are the $a=m-1$, $b=m+1$, 
cases of the following elements 
\ban 
E_{a\, b} &=& E_{a\, c} E_{c\, b} - q_c^{-1} E_{c\, b} E_{a\, c}, \\  
E_{b\, a} &=& E_{b\, c} E_{c\, a} -    q_c   E_{c\, a} E_{b\, c}, 
\ \ \ a<c<b. 
\nan 
The $\Zz$ grading of the algebra is specified such that 
the elements $K_a^{\pm 1}$, $\forall a\in {\bf I}$, 
and $E_{b\,  b+1}$, $E_{b+1\,  b}$,  $b\ne m$, are even, 
while $E_{m\, m+1}$ and $E_{m+1\, m}$ are odd.  
It is well known that $\Uq$ has the structure of a $\Zz$ graded 
Hopf algebra, with a co - multiplication 
\ban 
\Delta(E_{a\, a+1}) &=& E_{a\,  a+1} \otimes
K_a K_{a+1}^{-1} + 1 \otimes E_{a\, a+1}, \\ 
\Delta(E_{a+1\, a}) &=& E_{a+1\, a }\otimes 1 + K_a^{-1} K_{a+1}  
\otimes E_{a+1\, a}, \\
\Delta(K_a^{\pm 1}) &=&K_a^{\pm 1}\otimes K_a^{\pm 1},
\nan 
co - unit
\ban 
\epsilon(E_{a\, a+1})&=&\epsilon(E_{a+1\, a})=0, 
\ \ \forall a\in{\bf I}', \\ 
\epsilon(K_b^{\pm 1})&=&1,  \ \ \ \forall b\in{\bf I}, 
\nan  
and anti - pode
\ban 
S(E_{a\, a+1}) &=& - E_{a\, a+1} K_a^{-1} K_{a+1}, \\ 
S(E_{a+1\, a}) &=& - K_a K_{a+1}^{-1}E_{a+1\, a}, \\  
S(K_a^{\pm 1}) &=&K_a^{\mp 1}\otimes K_a^{\mp 1}. 
\nan

The quantum supergroup $\Uq$ has various $\Zz$ graded Hopf subalgebras, 
which are useful for analysing the representation theory.   
{}From equation (\ref{quantum}), we can easily see that 
{\em \begin{description} 
\item $\{ K_a \ | \ a\in{\bf I}\}$  generate 
an abelian Lie algebra 
\ban 
[K_a, \ K_b]&=& K_a K_b - K_b K_a =0.  
\nan 
\item $\{ K_a, \  a\in{\bf I}; \  E_{b\ b+1},\ E_{b+1\ b}, \ 
| \ m+n-1> b\in{\bf I}'\}$ 
generate a subalgebra $U_q(gl(m|n-1)\oplus gl(1))$. 
\item $\{ K_a \ | a\in{\bf I}\}$ $\cup$ $\{ E_{b\ b+1} \ | \ b\in{\bf I}'\}$  
generate a subalgebra $U_q({\frak b}_+)$, 
which is isomorphic to the quantized universal 
enveloping algebra of ${\frak b}_+$.   
\item Fix any subset ${\bf\Theta}\subset{\bf I}'$, then  
$\{K_a^{\pm 1}\ | \ a\in {\bf I}\}$ $\cup$ 
$\{ E_{b \ b+1} \ | \ b\in{\bf I}'\}$  $\cup$ 
$\{ E_{c+1\  c} \ | \ c\in{\bf\Theta}\}$ 
generate a ( ${\bf\Theta}$ dependent) subalgebra $U_q(\frak p)$, 
which is isomorphic to the quantized universal enveloping algebra of a 
standard parabolic subalgebra of $\gl$. 
\end{description}  }

Properties of the $E_{a b}$ were studied extensively in \cite{II}. 
We recall some of them below, which will be used repeatedly in the 
\BBW construction.  
\begin{lemma}\label{Eab} 
 1.  Assume $a<b$, then
\ban  
[ E_{a\, b},\ E_{c\, c+1} \}&=& 0, \\  
{[} E_{b\, a},\ E_{c+1\, c} \}&=&0, \ \ \ 
a\ne c,\ c+1,  \  \&  \  b\ne c,\ c+1, \\ 
{[} E_{a\, b},\ E_{c+1\, c} \}&=& \delta_{b\, c+1} 
E_{a\,c} K_c K_{c+1}^{-1} q_c^{-1} 
- \delta_{a\, c} (-1)^{\delta_{c\, m}} E_{c+1\, b} K_c^{-1} K_{c+1}, \\     
{[} E_{b\, a},\ E_{c\, c+1} \}&=& \delta_{a\, c} E_{b\, c+1} 
K_c K_{c+1}^{-1} q_{c+1} - \delta_{b\, c+1} (-1)^{\delta_{c\, m}} 
E_{c\, a} K_c^{-1} K_{c+1},\\    
& &a\ne c,\ c+1,\ \  \mbox{or}\ \   \ b\ne c,\ c+1. 
\nan 
2. 
\ban 
{[} E_{a\, b},\ E_{b\, a} \} &=& (K_a K_b^{-1} 
- K_a^{-1} K_b)/(q_a - q_a^{-1}), \\  
{[}E_{a\, c}, \  E_{c\, b}\}&=& 
     \left\{\begin{array}{l l}  
     E_{a\, b} K_c K_b^{-1} q_b, & a>b>c,\\ 
     E_{a\, b} K_c^{-1}  K_b,    & b>a>c, \\ 
     E_{a\, b} K_a^{-1} K_c,     & b<a<c,\\ 
     E_{a\, b} K_b K_c^{-1} q_b^{-1}, & a<b<c,  
     \end{array}\right.\\   
E_{c\, a} E_{c\, b}&=&
(-1)^{([a]+[c])([b]+[c])} q_c E_{c\, b} E_{c\, a},\\  
E_{b\, c} E_{a\, c}&=& 
(-1)^{([a]+[c])([b]+[c])} q_c^{-1} E_{a\, c} E_{b\, c},  
\ \ \ a<b<c, \ \mbox{or} \  \  b>a>c,  \\ 
{[}E_{c\, a},\ E_{c\, b}\}&=&[E_{a\, c}, \ E_{b\, c}\}=0,  
\ \ \ a<c<b, \ \mbox{or} \  \ a>b>c. 
\nan 
3.  Assume that $a<b$, $c<d$, and no two of $a, \ b, \ c$ and $d$ are 
equal.  Introduce the sets  $S(a,\, b)=\{a, a+1, ..., b\}$, and 
$S(c,\,  d)=\{c, c+1, ..., d\}$.   
\ban 
\mbox{If} \ \ \ S(a,\, b)\cap S(c,\,  d)
          &=& \emptyset , \ S(a,\, b),\ {\mbox or} \ \ S(c,\,  d),\\  
\mbox{then} \ \ \ [E_{a\, b},\ E_{c\, d}\}=[E_{a\, b},\ E_{d\, c}\}&=&
[E_{b\, a}, \ E_{c\, d}\}=[E_{b\, a}, \ E_{d\, c}\}=0. 
\nan   
\end{lemma}

Define 
\ban
\DE_{a b} &=& \DE_{a c} \DE_{c b} -    q_c   \DE_{c b} \DE_{a c}, \\
\DE_{b a} &=& \DE_{b c} \DE_{c a} - q_c^{-1} \DE_{c a} \DE_{b c}, \ \ \ a<c<b.
\nan
Then the $\DE_{a\,  b}$ are related to $E_{a\, b}$ with the help of 
the anti-pode 
\ban 
S(E_{a b}) &=& - \DE_{a b} K_a^{-1} K_b, \\ 
S(E_{b a}) &=& - K_a K_b^{-1} \DE_{b a},   \ \ \ a<b.
\nan 
Define the bilinear adjoint action $Ad:  \Uq\otimes \Uq\rightarrow \Uq$ by 
\ban 
Ad_x (y) &=& \sum_{(x)} (-1)^{[x_{(2)}][y]} x_{(1)} y S(x_{(2)}),  
\nan 
then $\Uq\rightarrow Ad_{\Uq}$,  $x\mapsto Ad_x$, yields an 
algebra homomorphism, where Sweedler's sigma notation is employed.  

Set 
\ban 
X_a &=& - \DE_{a\ m+n} K_a^{-1} K_{m+n},\\  
Y_a&=&  E_{m+n,  a},      \ \ \ a\in{\bf I}',  
\nan 
and denote by $\cal X$ the linear span of $\{ X_a\ | a\in{\bf I}'\}$, 
and by $\cal Y$ that of $\{ Y_a\ | a\in{\bf I}'\}$.  Then  
\begin{lemma}  
$\cal X$ and $\cal Y$  are closed  
under the adjoint action $Ad$ of the 
the subalgebra $U_q(gl(m|n-1)\oplus gl(1))$ $\subset$  $\Uq$,  
\ban
Ad_{K_a} X_b&=& q_a^{\delta_{a b}} X_b, \\ 
Ad_{e_c} X_b&=& \delta_{c+1, b} X_c, \\ 
Ad_{f_c} X_b&=& \delta_{c b}    X_{c+1},\\      
Ad_{K_a} Y_b&=& q_a^{-\delta_{a b}} Y_b, \\ 
Ad_{e_c} Y_b&=& - q_{c+1} \delta_{c b} Y_{c+1}, \\    
Ad_{f_c} Y_b&=& - q_{c+1}^{-1} \delta_{a+1, b} Y_c, \\  
Ad_{K_{m+n}}  X_b&=& q_{m+n}^{-1} X_b, \\ 
Ad_{K_{m+n}} Y_b &=& q_{m+n} Y_b,   \nan 
where $e_c=E_{c\, c+1}$, $f_c=E_{c+1\, c}$, $a, b, c\in{\bf I}'$, 
$\ c<m+n-1$. \end{lemma}  

Regarded as $U_q(gl(m|n-1)\oplus gl(1))$ modules, 
both  $\cal X$ and $\cal Y$ are irreducible, and  
are dual to each other in the following sense: 
the bilinear pairing ${\cal Y}\otimes{\cal X}\rightarrow$ 
${\bf\Lambda}_\infty$ defined by 
$(Y_a, \ X_b)=\delta_{a b}$ identifies $\cal Y$  with the dual 
vector space ${\cal X}^*$ of ${\cal X}$.  
This also defines a $U_q(gl(m|n-1)\oplus gl(1))$ 
module isomorphism 
\ban 
(Ad_u (Y_a), \ X_b) &=& (-1)^{[u]} (Y_a,  \ Ad_{S(u)}(X_b)), 
\ \ \ u\in U_q(gl(m|n-1)\oplus gl(1)). 
\nan  
Therefore,  
$\bar{\cal{C}}:= \sum_{a\in{\bf I}'}(-1)^{[b]+1} S^2(Y_a)\otimes X_a$, 
defines an invariant of the quantum subgroup 
\ban 
Ad_u(\bar{\cal{C}})&=& \sum_{a\in{\bf I}'} (-1)^{[u_{(2)}][a+1]+[b]+1}
                 Ad_{u_{(1)}}(S^2(Y_a))\otimes Ad_{u_{(2)}}(X_a)\\ 
     &=&\epsilon(u) \bar{\cal{C}},   \ \ \ u\in U_q(gl(m|n-1)\oplus gl(1)).   
\nan 
{}From this fact we can easily deduce that 
\begin{lemma}\label{commuting}
The ${\cal C}$ defined by  
\ban 
{\cal C}&=& \sum_{a\in{\bf I}'} (-1)^{[b]+1} Y_a\otimes S^{-1}(X_a),   
\nan 
satisfies 
\ba 
[\Delta'(u),  \ {\cal C}]&=&0,  
\ \ \ \forall u\in U_q(gl(m|n-1)\oplus gl(1)), 
\na 
where $\Delta'$ is the opposite co - multiplication. 
\end{lemma}  

The finite dimensional irreducible representations of $\Uq$ defined over 
the complex field were studied systematically in \cite{II}. The main 
conclusions still apply to the present situation.  We have 
{\em \begin{description}
\item Every finite dimensional irreducible $\Uq$ module admits a basis, 
      relative to which the $K_a$ are diagonal. 
\item Every finite dimensional irreducible $\Uq$ module is of highest weight 
      type and is uniquely ( up to isomorphisms ) characterized by a 
      highest weight. 
\item An irreducible $\Uq$ module $V(\lambda)$ with the maximal 
      vector $v^\lambda_+$ 
  \ban E_{a\ a+1}\ v^\lambda_+&=& 0, \ \ \ a\in{\bf I}',\\ 
       K_b\ v^\lambda_+&=& q_b^{\lambda_{b}} v^\lambda_+,  \ \ \ b\in{\bf I}. 
   \nan    
   is finite dimensional iff $\lambda$ satisfies 
   $\lambda_a - \lambda_{a+1}\in\Z_+$,  $a\ne m$. 
\item When $V(\lambda)$ is finite dimensional, it has the same weight 
      space decomposition as that of the corresponding irreducible 
      $\gl$ module with highest weight $\lambda$.   
\end{description} }

\subsection{A realization of $\Uq$ on projective superspace}
Before embarking on the \BBW construction for the quantum 
supergroup $\Uq$, we consider first a simple realization for 
it in terms of difference operators on the projective superspace 
$\C P^{m|n-1}$.  
   
Let $z$ be a variable living in $\C_c$. Define a difference operator 
$\nabla_z$ on analytic functions by 
\ban 
\nabla_z f(z) &=& { {f(q z) - f( q^{-1} z)}\over{ z ( q - q^{-1} )} }.  
\nan 
Then 
\ban 
\nabla_z(f(z) h(z))&=&\nabla_z f(z)\ q^{d_z} h(z) 
                    + q^{-d_z}f(z)\ \nabla_z h(z)\\ 
                   &=&\nabla_z f(z)\ q^{-d_z} h(z) 
                    + q^{d_z} f(z)\ \nabla_z h(z),  
\nan 
where $$d_z=z{ {d}\over{d z}}$$ is the scaling operator.  
For a Grassmannian variable $\theta$, we denote 
\ban 
\nabla_\theta&=&{{d}\over{d \theta}}, \\ 
d_{\theta}&=& \theta{ {d}\over{d \theta}}. 
\nan 

Consider $( W_a )_{a\in{\bf I}}= (\zeta_i; w_\mu)\in\C^{m|n}$. 
It is clear that  
\ban 
[\nabla_{w_\mu},\ \nabla_{w_\nu}]&=&0, \\  
{[}\nabla_{w_\mu},\ \nabla_{\zeta_i}]&=&0, \\
\{\nabla_{\zeta_i},\ \nabla_{\zeta_j}\}&=&0, \ \ \  \forall i, j, \mu, \nu.
\nan 

The quantum supergroup $\Uq$ can be realized in terms of the 
difference operators. Direct calculations can easily establish 
that the following operators satisfy the defining relations 
(\ref{quantum}) of $\Uq$
\ba 
{\cal E}_{a\ a+1} &=& - W_{a+1} \nabla_{W_a},\nonumber \\ 
{\cal E}_{a+1\ a} &=& -(-1)^{[a+1]([a]+1)}\, W_{a} \nabla_{W_{a+1}}, 
\ \ \ a\in{\bf I}', \nonumber \\  
{\cal K}^{\pm 1}_b &=& q^{\pm c} q_b^{\mp d_{W_b}}, 
                       \ \ \ b\in{\bf I},   
\na 
where $c\in\C_c$ is an arbitrary but fixed even number.  The homogeneous 
polynomials in the components of $( W_a )$ of a given degree 
furnishes an irreducible module over $\Uq$ in this realization. 
Such irreducible modules can all be obtained through repeatedly 
tensoring the contragredient vector module with itself. 

To turn this realization of $\Uq$ into one on $\C P^{m|n-1}$, we 
consider, e.g., the coordinate patch $U_{m+n}\subset\C P^{m|n-1}$. 
Let \ban 
( Z_a )_{a\in{\bf I}'} &=&  (\theta_i; z_\mu)_{1\le i\le m; \ m+1\le \mu\le m+n-1}, \\ 
\theta_i &=& \zeta_i/w_{m+n}, \ \ \ 1\le i\le m, \\ 
z_{\mu}&=& w_{\mu}/w_{m+n},   \ \ \ m+1\le \mu\le m+n-1.   
\nan 
Set $\nabla_{Z_a} = \nabla_a$.  Then 
\ban 
d_{W_a} &=& d_a ,  \\
W_a \nabla_{W_b} &=& Z_a \nabla_{b} ,  \\
W_{m+n} \nabla_{W_b}&=& \nabla_{b} , \\ 
W_a\nabla_{W_{m+n}}&=& Z_a W_{m+n}\nabla_{W_{m+n}}, 
\ \ \ \forall a, b \in{\bf I}',   
\nan 
where $W_{m+n}\nabla_{W_{m+n}}$ can be expressed as 
\ban 
W_{m+n}\nabla_{W_{m+n}}&=& { {q^{d_{W_{m+n}}} - 
q^{-d_{W_{m+n}}}}\over{q - q^{-1}} }.  
\nan  
Recall that $\prod_{a\in{\bf I}}(K_a)^{(-1)^{[a]}}$ is a central element 
of $\Uq$.  In an irrep realized by degree $k$ homogeneous polynomials 
in $( W_a )$, it takes the eigenvalue $q^{(m-n) c - k}$. Thus in this 
irrep, ${\cal K}_{m+n}$ and also $W_{m+n}\nabla_{W_{m+n}}$
 can be expressed in terms of 
the ${\cal K}_a$, $a\in{\bf I}'$.

Now the irrep can be realized solely in terms of the new variables 
$( Z_a )$.  We collect the results into 
\begin{proposition} 
The polynomials in  components of $( Z_a )_{a\in{\bf I}'}\in U_{m+n} 
\subset \C P^{m|n-1}$   
of degrees less or equal to $k$ furnishes an irreducible 
$\Uq$ module with the generators realized in terms of 
difference operators on $U_{m+n}$ by 
\ba
{\cal K}_a &=& q^c\, q_a^{- d_a},\ \ \ a\in{\bf I}',  \nonumber \\
{\cal E}_{a\ a+1} &=& - Z_{a+1} \nabla_a,\nonumber \\
{\cal E}_{a+1\ a} &=& -(-1)^{[a+1]([a]+1)}\, W_{a} \nabla_{a+1},
                       \ \ \ a+1\in{\bf I}', \nonumber \\ 
{\cal K}_{m+n}&=&q_{m+n}^{-k} q^c 
\prod_{a\in{\bf I}'}q_{m+n}^{d_a},    \label{subgroup}\\  
{\cal E}_{m+n-1\ m+n}&=& -\nabla_{m+n-1}, \nonumber \\
{\cal E}_{m+n\ m+n-1}&=& -  z_{m+n-1} 
{ {q^{-c} {\cal K}_{m+n} - q^c  {\cal K}_{m+n}^{-1} }
\over{q_{m+n} - q_{m+n}^{-1}} }.  
\na
Such a realization can be constructed for each coordinate patch of 
$\C P^{m|n-1}$, and the various realizations are related by 
coordinate changes on the projective superspace. 
\end{proposition}

The ${\cal E}_{a\ a+1}$,  ${\cal E}_{a+1\ a}$, $a<m+n-1$, 
and ${\cal K}_b$, $b\in{\bf I}$ given by equation (\ref{subgroup}) 
realize the $U_q(gl(m|n-1)\oplus gl(1))$ subalgebra of $\Uq$.  
We will denote this realization by $\Upsilon(U_q(gl(m|n-1)\oplus gl(1)))$, 
and introduce the algebra homomorphism 
 $\Upsilon: U_q(gl(m|n-1)\oplus gl(1))\rightarrow 
\Upsilon(U_q(gl(m|n-1)\oplus gl(1)))$ defined by 
$\Upsilon({E}_{a\ a+1})={\cal E}_{a\ a+1}$, etc.. 

We now construct the tensor operator ${\cal Y}$ of the 
subalgebra in this representation, the components of which will 
be denoted by $y_a$.  
Note that the complicated factor multiplying $z_{m+n-1}$       
in the expression of ${\cal E}_{m+n\ m+n-1}$ commutes 
with the subalgebra, thus we can ignore it without 
affecting the tensorial properties of the tensor operator. 
A simple calculation gives 
\ban 
y_a&=& (-1)^{[a]+1} Z_a q^{-\sum_{b=a+1}^{m+n-1} \{(-1)^{[b+1]} + d_b\}}, 
\ \ \ a\in{\bf I}'. 
\nan 
The corresponding $\cal C$ operator, which will be  denoted by 
$\cal O$, is given by 
\ba 
{\cal O}&=& \sum_{a\in{\bf I}'} (-1)^{[a]+1} y_a\otimes E_{a\, m+n}. 
\na
An immediate consequence of Lemma \ref{commuting} 
is that 
\begin{lemma}\label{Upsilon}   
\ba 
[(\Upsilon\otimes id) \Delta' (u),\ {\cal O}]&=&0, \ \ \ \forall u\in 
U_q(gl(m|n-1)\oplus gl(1)),  
\na 
where $\Delta'$ represents the opposite co - multiplication.  
\end{lemma}

Let us re - write $\cal O$ as  
\ban {\cal O}&=& {\cal O}' ( q^{-d_{m+n-1}}\otimes 1 ) 
+ z_{m+n-1} \otimes E_{m+n-1\ m+n}. \nan 
Then ${\cal O}'$ does not depend on $z_{m+n-1}$, and 
\ban 
{\cal O}' ( z_{m+n-1} \otimes E_{m+n-1\ m+n} ) 
&=& q^{-1}\, ( z_{m+n-1} \otimes E_{m+n-1\ m+n} ) {\cal O}'. 
\nan 
Introduce the formal power series  
\ba 
g(Z)&=& exp_q (\cal O ), \label{operator}\\ 
exp_q ( x ) &=&\sum_{k=0}^\infty {{x^k}\over{[k]_q!}}.\nonumber  
\na 
The $g(Z)$  is well behaved when acting on $\LC[[Z]]\otimes V$ 
if the $\Uq$ module $V$ is finite dimensional.  
Consider the action of  ${\cal O}^k$ on $1\otimes v$ $\in$ 
$\LC[[Z]]\otimes V$, where $v$ is an arbitrary element of $V$. 
An easy induction can establish that 
\ban 
&&\left[ {\cal O}' ( q^{-d_{m+n-1}}\otimes 1 ) 
+ z_{m+n-1} \otimes E_{m+n-1\ m+n}
\right]^k (1\otimes v)\\ 
&&=\sum_{l=0}^{k}{ {[k]_q!}\over{[k-l]_q!\, [l]_q!} } ({\cal O}')^{k-l}
(z_{m+n-1}^l\otimes E_{m+n-1\, m+n}^l v),  
\nan    
which in turn leads to  
\ban 
exp_q({\cal O}) \, (1\otimes v)&=&exp_q({\cal O}')\,  exp_q(z_{m+n-1} \otimes 
E_{m+n-1\, m+n} )\, (1\otimes v). 
\nan

The action of powers of ${\cal O}'$ on $z_{m+n-1}^k\otimes v$ can be 
similarly simplified. Such a process can be continued, and we 
eventually arrive at
\begin{lemma}\label{computation} 
\ba exp_q({\cal O}) \, (1\otimes v)&=&exp_q({\cal O}_1)
\, exp_q({\cal O}_2) ... exp_q({\cal O}_{m+n-1}) (1\otimes v), 
\ \ \ \forall v\in V,  \na  
where
\ban  
{\cal O}_a&=&\tilde{Z}_a\otimes E_{a\, m+n}, \\  
\tilde{Z}_a&=& q^{-\sum_{b=a+1}^{m+n-1} (-1)^{[b+1]}}\ Z_a, 
\ \ \  a\in{\bf I}'.     
\nan    
\end{lemma}  

\subsection{\BBW construction} 
With the preparations of the last two subsections, we can now 
carry out the \BBW construction for $\Uq$. 
Let $U_q(\frak{p})$ be the parabolic subalgebra of $\Uq$ generated by 
\ban 
K_a^{\pm 1}, \ \ a\in{\bf I}; & 
E_{b\, b+1}, \ E_{c+1\, c},&  
\ \ b, \, c \in{\bf I}', \ c< m+n-1,  
\nan 
subject to the appropriate relations of (\ref{quantum}).  
Then from results of \cite{II} we can deduce that 
\ban 
\Uq&=&{\cal U}^{(-)}\, U_q(\frak{p}), 
\nan 
where ${\cal U}^{(-)}$ is the subalgebra generated by 
$\{ E_{m+n\, a}\, |\, a\in{\bf I}'\}$.  

The induced module construction of subsection 2.2 can be generalized 
to the present case without much difficulty.  Given any irreducible 
$U_q(\frak{p})$  module $V_0$, we construct the {\em infinite 
dimensional} $\Uq$ module 
\ban 
{\bar V}&=&\Uq\otimes_{U_q(\frak{p})} V_0.
\nan 
When ${\bar V}$ is irreducible, we set $V={\bar V}$.  If it is not 
irreducible, then the given condition that $V_0$ is an 
irreducible $U_q(\frak{p})$ module leads to the existence of a unique 
maximal submodule $M$ of ${\bar V}$.  Now $V={\bar V}/M$ is irreducible. 

When $V_0$ is finite dimensional, it admits a unique maximal vector 
$v_+$ such that 
\ban E_{a \, a+1}\, v_+&=&0, \ \ \  a\in{\bf I}',\\  
     K_b \, v_+&=&q_b^{\lambda_b}\, v_+, \ \ \  b\in{\bf I}, 
\nan 
and the $\lambda_a$ must satisfy the condition that 
$\lambda_a - \lambda_{a+1}$ $\in$ $\Z_+$,  $m\ne a<m+n-1$. 
If the highest weight satisfies the further condition that 
$\lambda_{m+n-1} -\lambda_{m+n}$ $\in$  $\Z_+$, then the resultant 
irreducible $\Uq$ module $V$ is finite dimensional. In this 
way, we can construct all the finite dimensional 
irreducible $\Uq$ modules.  

We will denote by $V(\lambda)$ the finite dimensional irreducible 
$\Uq$ module with highest weight $\lambda$.  Decompose it into 
eigenspaces of $K_{m+n}$, 
\ban 
V(\lambda)&=&\oplus_{k=0}^L \ V^{(k)}(\lambda), \ \ \ L<\infty,  
\nan 
where 
\ban 
K_{m+n}\, v&=& q_{m+n}^{k+\lambda_{m+n}} \, v, 
\ \ \ \forall v\in V^{(k)}(\lambda). 
\nan   

We choose a basis for each subspace $V^{(k)}(\lambda)$. Then the union 
of all of them furnishes a basis $\{ w^\alpha\, |\, \alpha=1, 2, ..., 
dim V(\lambda)\}$ for $V(\lambda)$.  We order the $w^\alpha$ 
in such a way that $v^i=w^i$, $i=1, 2, ..., dim V^{(0)}(\lambda)$, 
are the basis elements of $V^{(0)}(\lambda)$.  
Let $\{ \langle w^\alpha \, | \  \rangle\}$ be the dual basis of 
$V(\lambda)^*$, i.e., 
$\langle w^\alpha  | w^\beta\rangle =\delta_{\alpha \, \beta}$.   

Similar to the classical case, we define  
a linear map $\Xi:$  $V(\lambda)$
$\rightarrow$   $\LC[[Z]]_L\otimes V_0(\lambda)$ by
\ba 
\Xi_w(Z)&=& \sum_{i=1}^{d_0} (-1)^{[v^i](1+[w])} \langle v^i |  
exp_q({\cal O}) \, ( 1 \otimes w) \rangle \otimes v^i, \\ 
& &d_0=dim V^{(0)}(\lambda),\nonumber   
\na
where $exp_q(\cal O)$ is given in (\ref{operator}). 
Since $V(\lambda)$ is finite dimensional, 
$exp_q({\cal O})(1\otimes w)$ 
is well defined for all $w\in V(\lambda)$. 

Using Lemma 3, we can easily see that the elements of
\ban 
\{\Xi_\alpha(Z)=\Xi_{w^\alpha}(Z)\, | \alpha=1, 2, ...,dim V(\lambda)\} 
\nan 
are linearly independent. Their linear span  
$V^\lambda_Z$ $\subset$ $\LC[[Z]]_L\otimes V_0(\lambda)$ yields 
a $\Uq$ module with the module action defined by  
\ban
u\circ \Xi_w(Z) &=& \Xi_{u w}(Z), \ \ \   \forall u\in\Uq,  
\nan
which is isomorphic to $V(\lambda)$ itself. 

The quantum supergroup  $\Uq$ can now be realized in terms of 
operators acting on $V^\lambda_Z$. We have  
\begin{proposition} 
Let $\pi$ denote the irreducible representation of $\Uq$ afforded 
by  the irreducible module $V^\lambda_Z$. Then  
\ban  
\pi(E_{a\,  a+1})&=& -Z_{a+1} \nabla_a\otimes 1 + 
   q_a^{-d_a} q_{a+1}^{d_{a+1}}\otimes \pi_0(E_{a\,  a+1}), \\
\pi(E_{a+1\, a})&=& -Z_{a} \nabla_{a+1}\otimes q_a^{d_a} 
q_{a+1}^{-d_{a+1}}    
+ 1\otimes \pi_0(E_{a+1\, a}), \ \ \ a<m+n-1, \\ 
\pi(K_b)&=& q_b^{-d_b}\otimes \pi_0(K_b), \ \ \ b\le m+n-1,\\ 
\pi(K_{m+n})&=& q_{m+n}^{\sum_{a=1}^{m+n-1} d_a} \otimes \pi_0(K_{m+n}), \\ 
\pi(E_{m+n-1\, m+n})&=& \nabla_{m+n-1}\otimes 1, \\ 
\pi(E_{m+n\, m+n-1})&=& 
{ {1}\over{q_{m+n-1} - q_{m+n-1}^{-1}} } \left\{ 
 (z_{m+n-1}\otimes \pi_0(K_{m+n-1}))\pi(K^{-1}_{m+n})\right.\\ 
&-&\left. (z_{m+n-1}\otimes \pi_0(K^{-1}_{m+n-1}))\pi(K_{m+n})\right\}\\
&+& \sum_{a=1}^{m+n-2} {\tilde Z}_a q_{m+n}^{-1-\sum_{b=1}^a d_b} 
\otimes \pi_0(E_{a\, m+n-1}\, K_{m+n-1} K_{m+n}^{-1}),  
\nan      
where $\pi_0$ denotes the irreducible representation of
$U_q(\frak{p})$ afforded by $V^{(0)}(\lambda)$. 
\end{proposition} 
{\em Proof}:   We consider the $U_q(gl(m|n-1)\oplus gl(1))$ 
subalgebra first. Observe that 
\ban  
(\Upsilon\otimes id) \Delta'  (E_{a\, a+1})\, (1 \otimes w) 
&=&  1\otimes E_{a\, a+1}\, w, \\ 
(\Upsilon \otimes id)\Delta' (E_{a+1\, a})\, (1\otimes w)  
&=&  1\otimes E_{a+1\, a}\, w,
\ \ \ \forall w\in V(\lambda),  \ \    a<m+n-1.  
\nan 
By applying Lemma \ref{Upsilon}, we obtain 
\ban 
E_{a\, a+1} \circ \Xi_w(Z)&=&  
\sum_{i=1}^{d_0} (-1)^{[v^i](1+\delta_{a m}+ [w])} \langle v^i |  
exp_q(\cal O)\, ( 1 \otimes E_{a\, a+1}\,  w) \rangle \otimes v^i\\ 
&=& \sum_{i=1}^{d_0} \langle v^i |({\cal E}_{a\, a+1}\otimes 1 
+ {\cal K}_a {\cal K}_{a+1}^{-1}\otimes E_{a\, a+1})\, exp_q({\cal O})\, 
( 1 \otimes w) \rangle \\ 
&\otimes &v^i \, (-1)^{[v^i](1+\delta_{a m}+ [w])}\\ 
&=& ({\cal E}_{a\, a+1}\otimes 1  
+ {\cal K}_a {\cal K}_{a+1}^{-1}\otimes \pi_0(E_{a\, a+1}) )\, \Xi_w(Z),  
\nan
which immediately yields the expression for $\pi(E_{a\, a+1})$. 
In exactly the same way, we can obtain the realizations 
of the other generators of this subalgebra. 

The realization of $E_{m+n-1\, m+n}$ can be easily obtained 
from Lemma \ref{computation}. To prove the formula for 
$E_{m+n\, m+n-1}$, however, considerable effort is required. 
We will need the following technical results  
\ba 
& &[exp_q({\cal O}_{m+n-1}), \ 1\otimes E_{m+n\, m+n-1}]\nonumber\\
&=&{ {1}\over{q_{m+n} - q_{m+n}^{-1}}}   \left[  
(z_{m+n-1}\otimes K_{m+n-1} ) exp_q({\cal O}_{m+n-1}) 
 (1\otimes K_{m+n}^{-1})\right. \nonumber \\ 
&-&\left.  (z_{m+n-1}\otimes K^{-1}_{m+n-1} ) exp_q({\cal O}_{m+n-1}) 
  ( 1\otimes K_{m+n}) \right],\label{Omn}  \\  
& &{[} exp_q({\cal O}_a), \ 1\otimes E_{m+n\, m+n-1}]\nonumber\\ 
&=& ({\tilde Z}_a q^{-1}_{m+n}\otimes E_{a\, m+n-1}) 
exp_q({\cal O}_a) (1\otimes K_{m+n-1} K_{m+n}^{-1}), \label{Oa}    
\na 
which can be proven by repeatedly applying Lemma \ref{Eab}. 

Set 
\ban 
\eta&=& \sum_{i=1}^{d_0} \langle v^i |  exp_q({\cal O}_1) ... 
exp_q({\cal O}_{m+n-2}) \left[exp_q({\cal O}_{m+ n-1}), 
\   1 \otimes E_{m+n\, m+n-1}\right](1 \otimes \,  w)
\rangle\\ 
&\otimes& v^i (-1)^{[v^i](1+[w])}\\
\zeta&=&\sum_{i=1}^{d_0} \langle v^i |\left[exp_q({\cal O}_1) ... 
exp_q({\cal O}_{m+ n-2}), \ 1 \otimes E_{m+n\, m+n-1}\right] 
exp_q({\cal O}_{m+ n-1}) (1 \otimes \,  w) \rangle\\ 
&\otimes& v^i\ (-1)^{[v^i](1 + [w])}. 
\nan  
Then 
\ban E_{m+n\, m+n-1}\circ\Xi_w&=&\eta + \zeta. \nan 
By using (\ref{Omn}), we can re - write $\eta$ as 
\ban 
\eta&=&
{ {(z_{m+n-1}\otimes \pi_0(K_{m+n-1})) 
  \pi(K^{-1}_{m+n}) 
-(z_{m+n-1}\otimes \pi_0(K^{-1}_{m+n-1}))
  \pi(K_{m+n})}\over{q_{m+n-1} - q_{m+n-1}^{-1}}} \Xi_w. 
\nan  
Lemma \ref{Eab} asserts that $E_{a\, m+n-1}$ (anti)commutes with 
all $E_{c\, m+n}$ if $c<a$.  Thus by using (\ref{Oa}) we 
arrive at 
\ban 
\zeta&=&\sum_{i=1}^{d_0}\sum_{a=1}^{m+n-2}  
\langle v^i |({\tilde Z}_a q^{-1}_{m+n}\otimes E_{a\, m+n-1}\, K_{m+n-1}) 
exp_q({\cal O}_1) ...exp_q({\cal O}_a)\\  
&\times& (1\otimes K_{m+n}^{-1})
exp_q({\cal O}_{a+1})... exp_q({\cal O}_{m+n-1}) 
(1 \otimes \,  w) \rangle \otimes v^i \ (-1)^{[v^i](1+[w])}\\ 
&=& \sum_{i=1}^{d_0}\sum_{a=1}^{m+n-2} 
\langle v^i |({\tilde Z}_a q_{m+n}^{-1-\sum_{b=1}^{m+n-1} d_b}
\otimes E_{a\, m+n-1}\, K_{m+n-1} K_{m+n}^{-1}) 
exp_q({\cal O}_1) ...exp_q({\cal O}_a)\\  
&\times& (q_{m+n}^{\sum_{b=1}^{m+n-1} d_b}\otimes 1)  
exp_q({\cal O}_{a+1})... exp_q({\cal O}_{m+n-1})
(1 \otimes \,  w) \rangle \otimes v^i \, (-1)^{[v^i](1+[w])} \\ 
&=& \sum_{a=1}^{m+n-2} 
({\tilde Z}_a q_{m+n}^{-1-\sum_{b=1}^{a} d_b}
\otimes \pi_0(E_{a\, m+n-1}\, K_{m+n-1}\, K_{m+n}^{-1} ))\,  \Xi_w.  
\nan
Adding the final expressions of $\eta$ and $\zeta$ together gives 
the desired realization of $E_{m+n\, m+n-1}$. 

\section{\normalsize CONCLUSION} 
As we have explained, the construction of the vector 
coherent states for the supergroup $GL(m|n)$ at 
the classical level can be regarded as one manifestation of 
a supersymmetric generalization of the celebrated \BBW theorem, 
and the vector coherent states themselves may be  interpreted 
as holomorphic sections of homogeneous supervector bundles. 
It should be possible to give the quantum coherent states   
a similar interpretation  within a yet to be fully developed 
framework of noncommutative geometry.  
The results on the quantum \BBW construction should  
feedback concrete useful information,  improving this 
framework itself. 
In a forthcoming  publication, we will develop 
a global version of the quantum \BBW  construction, the connection of 
which with the results presented here will also be clarified. 
In doing so,  results of references 
\cite{Schneider} \cite{Parshall} \cite{Noumi}\cite{Koorwinder} 
on quantum fiber bundles will come to play.

\vspace{2cm}

\noindent 
ACKNOWLEDGEMENTS:\\ 
This paper originated from some unpublished work of Professor A. J.
Bracken's, who also carefully read the first draft of the paper, 
and gave many instructive comments.

%\pagebreak 
 

\vspace{2cm}
\begin{thebibliography}{9999}
\bibitem{Bartocci} C. Bartocci, U. Bruzzo and D. Hernandez - Ruiperez, 
             {\em The geometry of supermanifolds}, Kluwer Academic Publishers, 
             London (1991). 
\bibitem{Bracken} A. J. Bracken,  unpublished notes (1989).   
\bibitem{Lohe} L. C. Biedenharn and M. Lohe, {\em An extension of 
       the Borel - Weil construction to the quantum group $U_q(n)$}, 
       Commun. Math. Phys. {\bf 146} (1992) 483.     
\bibitem{Biedenharn}L. C. Biedenharn and M. Lohe,
        {\em Quantum group symmetry and $q$ - tensor operator algebras}, 
        World Scientific, Singapore (1995).  
\bibitem{Bott} R. Bott, {\em Homogeneous vector bundles}, 
         Ann. Math. {\bf 66} (1957) 203. 
\bibitem{Cornwell} J. F. Cornwell, {\em Abstract carrier space 
        formalism for the irreducible tensor operators of compact 
        quantum group algebras}, St. Andrews Preprint. 
\bibitem{DeWitt} B. DeWitt, {\em Supermanifolds}, Cambridge Press, 
        London (1984). 
\bibitem{Jarvis} P. H. Dundi and P. D. Jarvis, 
	{\em Diagram and superfield techniques in the classical 
	 superalgebras},  J. Phys. {\bf A 14} (1981) 547.   
\bibitem{Koorwinder} M. S. Dijkhuizen and T. H. Koorwinder, 
         {\em Quantum homogeneous spaces, duality and quantum 2 - spheres}, 
         Geom. Dedicata {\bf 52} (1994) 291. 
\bibitem{Hecht} P. Exner, M. Harlicek and W. Lassner, 
		Czech. J. Phys., {\bf B26} (1976) 1213;\\ 
	 K. T. Hecht, {\em The vector coherent state method and 
         its applications to problems in higher symmetries}, 
        Lecture Notes in Physics, 
        {\bf 290} Springer - Verlag, Berlin(1987), and references
	therein.      
\bibitem{Kac} V. G. Kac, {\em Lie superalgebras}, 
              Adv. Math. {\bf 26} (1977) 8.   
\bibitem{Rowe} R. Le Blanc and D. J. Rowe, {\em Highest weight
               representations of $gl(m/n)$ and $gl(m+n)$}, 
	       J. Math. Phys., {\bf 30} (1989) 1415. 
\bibitem{Manin} Y. I. Manin,  {\em Gauge field theory and complex geometry}, 
                Springer - Verlag, Berlin (1988). 
\bibitem{Noumi} M. Noumi, H. Yamada and K. Mimachi, {\em Finite - dimensional 
           representations of the quantum group $GL_q(n,\, \C)$ and 
           zonal spherical functions on $U_q(n-1)\backslash U_q(n)$}, 
           Japanese J. Math., {\bf 19} (1993) 31. 
\bibitem{Parshall} B. Parshall and J. P. Wang, {\em Quantum linear groups},  
		  Memoirs Amer. Math. Soc., {\bf 89} No. 439 (1991) 1 - 157. 
\bibitem{Penkov} I. B. Penkov,  {\em Borel - Weil - Bott Theory for Classical 
         Lie Supergroups},  Sovr. Probl. Math. {\bf 32} VINITI, Moscow, 
         (1988) 71 - 124.   
\bibitem{Serganova} I. B. Penkov and V. Serganova, {\em Cohomology of 
        $G/P$ for Classical Complex Lie Groups of $G$ and Characters of 
        Some Typical $G$ - Modules},  Annales de L'institut Fourier, 
        {\bf 39} (1989) 846 - 873.  
\bibitem{Scheunert}M. Scheunert, {\em The theory of Lie superalgebras}, 
        Lecture Notes in Math., {\bf 716} Springer - Verlag, Berlin (1979). 
\bibitem{Schneider} H. J. Schneider, 
        {\em Principle homogeneous spaces for arbitrary Hopf algebras}, 
        Israel J. Math. {\bf 72} (1990) 196. 
\bibitem{II} R. B. Zhang, {\em Finite Dimensional Representations of 
         the Quantum Supergroup $U_{q}(gl(m/n))$}, 
        J. Math.  Phys., {\bf 34}(1993) 1236 -- 1254.
\end{thebibliography}
\end{document}